\documentclass[preprint,12pt]{elsarticle}
\usepackage{amsmath}
\usepackage{latexsym}
\usepackage{epsfig}
\usepackage{psfrag}
\usepackage{enumerate}
\usepackage[utf8]{inputenc}
\usepackage{newtxtext}
\usepackage{xcolor}
\usepackage[bb=dsserif]{mathalpha}
\usepackage{bm}

\makeatletter
\def\ps@pprintTitle{%
 \let\@oddhead\@empty
 \let\@evenhead\@empty
 \def\@oddfoot{\centerline{\thepage}}%
 \let\@evenfoot\@oddfoot}
\makeatother

\newtheorem{thm}{Theorem}[section]
\newtheorem{defin}{Definition}[section]
\newtheorem{lemma}{Lemma}[section]

\newtheorem{rem}{Remark}[section]

\newtheorem{prop}{Proposition}[section]

\def\sgn{{\rm sgn}}

\usepackage{amssymb}

\journal{Journal of Functional Analysis}

\begin{document}

\begin{frontmatter}



\title{On the value of the fifth maximal projection constant}

 \author[label1]{Beata~Der\c{e}gowska\footnote{B.D. is partially supported by National Science Center (NCN) grant no. 2021/05/X/ST1/01212. For the purpose of Open Access, the author has applied a CC-BY public copyright licence to any Author Accepted Manuscript (AAM) version arising from this submission.}}
 \author[label2]{Matthew~Fickus\footnote{The views expressed in this article are those of the authors and do not reflect the official policy or position of the United States Air Force, Department of Defense, or the U.S. Government.}}
 \author[label3]{Simon~Foucart\footnote{S.F. is partially supported by grants from the NSF (DMS-2053172) and from the ONR (N00014-20-1-2787).}}
 \author[label4]{Barbara~Lewandowska}
 
\address[label1]{Institute of Mathematics\\
Pedagogical University of Krakow, Podchorazych~2, Krakow, 30-084, Poland}
\address[label2]{Department of Mathematics and Statistics\\
Air Force Institute of Technology, Wright-Patterson AFB, OH 45433, USA}
\address[label3]{Department of Mathematics, Texas A\&M and Institute of Data Science\\
Texas A\&M University, College Station, TX 77843, USA}
\address[label4]{Faculty of Mathematics and Computer Science\\
Jagiellonian University, Lojasiewicza~6, Krakow, 30-048, Poland}

\begin{abstract}
Let $\lambda(m)$ denote the maximal absolute projection constant over real $m$-dimensional subspaces.
This quantity is extremely hard to determine exactly,
as testified by the fact that the only known value of $\lambda(m)$ for $m>1$  is $\lambda(2)=4/3$.
There is also numerical evidence indicating that $\lambda(3)=(1+\sqrt{5})/2$.
In this paper,
relying on a new construction of certain mutually unbiased equiangular tight frames,
we show that $\lambda(5)\geq 5(11+6\sqrt{5})/59 \approx 2.06919$.
This value coincides with the numerical estimation of $\lambda(5)$ obtained by B. L. Chalmers,
thus reinforcing the belief that this is the exact value of $\lambda(5)$.
\end{abstract}




\end{frontmatter}




\section{Introduction}
Let $X$ be a real Banach space and $Y\subset X$ be a finite-dimensional subspace.
Let $\mathcal{P}(X,Y)$ denote the set of all linear and continuous projections from $X$ onto $Y$,
recalling that an operator $P \colon X \rightarrow Y$ is called a \textit{projection} onto $Y$ if $P|_Y={\rm Id}_Y.$
We define the \textit{relative projection constant} of $Y$ by
\begin{equation*}
\lambda(Y,X) :=\inf\lbrace\|P\|:\;P\in\mathcal{P}(X,Y)\rbrace
\end{equation*}
and the \textit{absolute projection constant} of $Y$ by
\begin{equation}
\label{DefMAPC}
\lambda(Y) :=\sup\lbrace\lambda(Y,X):Y\subset X\rbrace.
\end{equation}
The literature also deals with the  \textit{maximal absolute projection  constant}, which is defined by
\begin{equation*}
\lambda(m) :=\sup \lbrace\lambda(Y):\; \dim(Y)=m \rbrace.
\end{equation*}
By the Kadec--Snobar theorem (see \cite{KS}),
we have $\lambda(m)\leq \sqrt{m}$.
Moreover, it has been shown in \cite{K} that this estimate is asymptotically the best possible.
However, the determination of the constant $\lambda(m)$ seems to be difficult: apart from $\lambda(1)=1$, the only known value of $\lambda(m)$ is $\lambda(2)=4/3$ --- this is \textit{Gr\"unbaum conjecture}, formulated in \cite{G} and proved in \cite{CL}.
Numerical computations presented in \cite{FS} indicate that $\lambda(3)$ should  equal $(1+\sqrt{5})/2$ --- this was stated,  with an erroneous proof, in \cite{KT}.
Other numerical experiments conducted by B. L. Chalmers (and unfortunately unpublished)  suggest that $\lambda(5)\approx 2.06919$.
In this article, we show that
$$
\lambda(5)\geq 5(11+6\sqrt{5})/59\approx 2.06919.
$$
Viewed in isolation, this could seem anecdotal.
However, several sources of evidence hint that this is the actual value of $\lambda(5)$.
This comes as a surprise, because it was growingly believed
that obtaining exact formulas for $\lambda(m)$ was an unreasonable quest.
Now there is hope that this quest could be realized after all.

To establish the announced lower bound,
we make a detour via maximal relative projection constants.
Recent results concerning maximal relative and absolute projection constants
can be found in \cite{B,B2,CLM,FS,S}.
Here, we only give the definition of the \textit{maximal relative projection constant} for $n \ge m$ as
\begin{equation*}
\lambda(m,n):=\sup\lbrace \lambda(Y, l_\infty^{(n)}):\; \dim(Y)=m \textrm{ and } Y\subset l_\infty^{(n)}\rbrace.
\end{equation*}
This is motivated by the fact that,
in the expression \eqref{DefMAPC} of $\lambda(m)$,
it suffices to take the supremum over finite-dimensional $l_\infty$ superspaces (see e.g. \cite[III.B.5]{W}),
so that  the nondecreasing sequence $( \lambda(m, n) )_{n \ge m}$ converges to $\lambda(m)$.
In reality,  there even is an $N\in\mathbb{N}$ such that $\lambda(m,n) = \lambda(m)$ for all $n\geq N$ (see \cite[Theorem~1.4]{B}).
Our estimation of $\lambda(m,n)$
will rely on the following result proved in \cite{CLe}.

\begin{thm}\label{lammbda}
For integers $n \ge m$, one has
\begin{equation*}
    \lambda(m,n)=\max\bigg\lbrace \sum_{i,j=1}^nt_it_j|U^\top U|_{ij}:t\in\mathbb{R}^n,\;\|t\|_2=1,U\in \mathbb{R}^{m\times n},\; UU^\top={\rm I}_m \bigg\rbrace .
\end{equation*}
\end{thm}

\noindent
Although this theorem provides an essential tool for estimating the maximal relative projection constants,
 computing their exact values remains a challenging problem,
carried out in just a few cases (see e.g. \cite{B,CLe,FS}).
One particular situation where an explicit formula is available
involves equiangular tight frames.
Let us recall that a system of unit (i.e., $l_2$-normalized) vectors $(v_1,\dots, v_n)$ in $\mathbb{R}^m$ is called {\it equiangular} if there is a constant $c\geq 0$ such that
\begin{equation*}
|\langle v_i,v_j \rangle|=c
\qquad \textrm{ for all } i,j\in \{1,\dots, n\},\;i\neq j.
\end{equation*}
It is called a \textit{tight frame} if
\begin{equation*}
VV^\top=\frac{n}{m}{\rm I}_m,
\end{equation*}
where $V$ is the matrix with columns $v_1,\dots, v_n$.
The system $(v_1,\dots, v_n)$ of unit vectors is called an {\it equiangular tight frame} if it is both  equiangular and a tight frame.
For an equiangular tight frame of $n$ unit vectors in $\mathbb{R}^m$,
it is well known (see e.g. \cite[Theorem 5.7]{FR}) that
\begin{equation*}
|\langle v_i,v_j \rangle|=\sqrt{\frac{n-m}{m(n-1)}}
\qquad \textrm{ for all } i,j\in \{1,\dots, n\},\;i\neq j.
\end{equation*}
The above-mentioned explicit formula is presented as part of the result below.
Built from Theorems 1 and 2 of \cite{KLL},
it appeared in a slightly different form as Theorem 5 in \cite{FS}.
A new self-contained proof is included later as an appendix.

\begin{thm}
\label{ETFs}
For integers $n\geq m$, the maximal relative projection constant $\lambda(m,n)$ is upper bounded by
$$
\delta_{m,n} := \frac{m}{n} \left( 1 + \sqrt{\frac{(n-1)(n-m)}{m}} \right).
$$
Moreover,  the equality $\lambda(m,n) = \delta_{m,n}$ occurs if and only if there is an equiangular tight frame for $\mathbb{R}^m$ consisting of $n$ unit vectors.
\end{thm}

\begin{rem}
\label{RemGam}
{\rm
We note in passing that $\delta_{m,n} < \sqrt{m}$ for $n \ge m > 1$
(thus providing another justification for Kadec--Snobar estimate).
This is seen by applying Cauchy--Schwarz inequality for the noncolinear vectors $[1, \sqrt{n-1}]$ and $[1, \sqrt{(n-m)/m} ]$ in
\begin{align*}
\delta_{m,n} & = \frac{m}{n} \bigg( 1 + \sqrt{n-1} \sqrt{\frac{n-m}{m}}  \bigg)
<  \frac{m}{n} \sqrt{1+n-1}  \sqrt{ 1+ \frac{n-m}{m} }
= \sqrt{m}.
\end{align*}
}
\end{rem}

In the rest of this paper,
we present new explicit lower bounds for $\lambda(m,n)$
under the condition that certain mutually unbiased equiangular tight frames for~$\mathbb{R}^m$ exist (see Theorem \ref{amazing}).
We then provide a construction of an infinite family of such mutually unbiased equiangular tight frames (see Theorem \ref{muetfFamily}).
Finally, combining these two ingredients,
we highlight the resulting estimation of $\lambda(5,16)$ to arrive at the promised lower bound for $\lambda(5)$,
conjectured to be its true value.

\section{The Lower Bound}

Before stating the main result,
we start with an observation about mutually unbiased equiangular tight frames,
formally defined below.

\begin{defin}\label{mu}
Two equiangular tight frames $(v_1,\dots,v_k)$ and $(w_1,\dots,w_l)$ for~$\mathbb{R}^m$
are mutually unbiased if there exists $c\in\mathbb{R}$ such that
\begin{equation*}
|\langle v_i,w_j \rangle|=c
\qquad \mbox{for all } i \in \{1,\dots,k\} \mbox{ and } j\in \{1,\dots,l\}.
\end{equation*}
\end{defin}

This definition generalizes a concept introduced in \cite{FM} so as to permit the case $k \neq l$.
We point out that the scalar $c$ is uniquely determined,
as also noted in~\cite{CpGaG}.

\begin{lemma}\label{beta}
The constant $c$ appearing in the definition of mutually unbiased equiangular tight frames for $\mathbb{R}^n$ necessarily  satisfies
$$
c=\frac{1}{\sqrt{m}}.
$$
\end{lemma}

\noindent
{\sc Proof.}
Let $(v_1,\dots,v_k)$ and $(w_1,\dots,w_l)$ be mutually unbiased equiangular tight frames for~$\mathbb{R}^m$
and let $V \in \mathbb{R}^{m \times k}$ be the matrix with columns $v_1,\dots,v_k$.
For any $j\in\{1,\dots,l\}$,
because the two frames are mutually unbiased, we have
\begin{equation*}
\|V^\top w_j\|_2^2
=\sum_{i=1}^k |\langle v_i,w_j \rangle|^2
=\sum_{i=1}^k c^2
=kc^2.
\end{equation*}
Since $(v_1,\dots,v_k)$ is a tight frame for~$\mathbb{R}^m$, we also have $VV^\top = (k/m){\rm I}_m$, and so
\begin{equation*}
\|V^\top w_j\|_2^2
=\langle V^\top w_j, V^{\top}w_j \rangle
=\langle w_j, VV^{\top}w_j \rangle
=\Big\langle w_j, \dfrac{k}{m}w_j \Big\rangle
=\dfrac{k}{m}\| w_j\|_2^2
=\dfrac{k}{m}.
\end{equation*}
It follows that $kc^2=k/m$, and hence $c=1/{\sqrt{m}}$, as claimed.\qed
\vskip .1in

We now present the main theorem of this section,
whose statement involves the quantity $\delta_{m,n}$ introduced in Theorem \ref{ETFs}.

\begin{thm}\label{amazing}
If mutually unbiased equiangular tight frames $(v_1,\dots,v_k)$ and $(w_1,\dots,w_l)$ for~$\mathbb{R}^m$ exist,
then the maximal relative projection constant $\lambda(m,k+l)$ is bounded below as
$$
\lambda(m,k+l) \geq  \frac{m-\delta_{m,k}\delta_{m,l}}{2\sqrt{m}-\delta_{m,k}-\delta_{m,l}}.
$$
\end{thm}

\noindent
{\sc Proof. }
Let $V \in \mathbb{R}^{m \times k}$ be the matrix with columns $v_1,\dots,v_k$ and $W \in \mathbb{R}^{m \times l}$ the matrix with columns $w_1,\dots,w_l$.
For any $\theta\in [0,\pi/2]$, let us consider the vector $t_\theta\in \mathbb{R}^{k+l}$
and the matrix $U_\theta  \in \mathbb{R}^{m \times (k+l)}$ defined, in block notation, by
\begin{equation}
\label{DeftU}
 t_\theta:= \begin{bmatrix} \cos\theta \dfrac{1}{\sqrt{k}} \mathbb{1}_k \\
 \hline \sin \theta  \dfrac{1}{\sqrt{l}}\mathbb{1}_l \end{bmatrix}
\qquad \mbox{and} \qquad
U_\theta:= \begin{bmatrix} \; \cos\theta \sqrt{\dfrac{m}{k}} V & \vline & \sin \theta \sqrt{\dfrac{m}{l}} W  \; \end{bmatrix},
\end{equation}
where $\mathbb{1}_n$ denotes the $n$-dimensional vector with all entries equal to $1$.
We observe that  $\|t_\theta\|_2=1$, that
\begin{equation*}
    U_\theta {U_\theta}^\top=\cos^2\theta \frac{m}{k}VV^\top + \sin^2 \theta \frac{m}{k} W W^\top=\cos^2\theta \, {\rm I}_m+\sin^2\theta \, {\rm I}_m = {\rm I}_m,
\end{equation*}
and that
\begin{equation}
\label{UTU}
  {U_\theta}^\top  U_\theta
  = \begin{bmatrix}
\cos^2 \theta \dfrac{m}{k} V^\top V &  \vline & \cos \theta \sin \theta \dfrac{m}{\sqrt{k l}} V^\top W\\
\hline
  \cos \theta \sin \theta \dfrac{m}{\sqrt{k l}} W^\top V & \vline & \sin ^2 \theta \dfrac{m}{l} W^\top W
\end{bmatrix}.
\end{equation}
Therefore, according to the expression of $\lambda(m,n)$ from Theorem \ref{lammbda},
we can make use of the tight frame and unbiasedness properties of $U$ and $V$
to obtain, with the shorthand notation $\phi_{m,n} := \sqrt{(n-m)/(m(n-1))}$,

\begin{align*}
    \lambda(m,k+l)& \geq \sum_{i,j=1}^{k+l} (t_\theta)_i (t_\theta)_j | {U_\theta}^\top  U_\theta|_{i,j} \\
    & =
    \cos^2 \theta  \frac{1}{k}  \times \cos^2\theta \frac{m}{k} \times k +
    \cos^2 \theta  \frac{1}{k} \times \cos^2\theta \frac{m}{k} \phi_{m,k} \times  k(k-1)\\
    & + \sin^2 \theta  \frac{1}{l}  \times \sin^2\theta \frac{m}{l} \times l +
    \sin^2 \theta  \frac{1}{l} \times \sin^2\theta \frac{m}{l} \phi_{m,l} \times  l(l-1)\\
   & +  2 \times \cos \theta \sin \theta \frac{1}{\sqrt{k l}} \times
    \cos \theta \sin \theta \frac{m}{\sqrt{k l}}  \frac{1}{\sqrt{m}}  \times k l \\
 & = \cos^4 \theta \bigg( \frac{m}{k} + \frac{m}{k} (k-1) \phi_{m,k} \bigg)+
     \sin^4 \theta \bigg( \frac{m}{l} + \frac{m}{l} (l-1) \phi_{m,l} \bigg)\\
  & + 2 \cos^2 \theta \sin^2 \theta \sqrt{m}\\
  & = \bigg(\frac{1+\cos(2\theta)}{2}\bigg)^2  \delta_{m,k}
  + \bigg(\frac{1-\cos(2\theta)}{2}\bigg)^2 \delta_{m,l}
  + \big( \sin(2\theta) \big)^2 \frac{\sqrt{m}}{2}.
\end{align*}
Since this is valid for any $\theta \in [0,\pi/2]$,
after setting $x:=\cos(2\theta)$,
we arrive at
\begin{align*}
\nonumber
  \lambda(m,k+l)&\geq \max_{x\in[-1,1]}\left( \frac{\delta_{m,k}(1+2x+x^2)}{4}+ \frac{\delta_{m,l}(1-2x+x^2)}{4}+\frac{\sqrt{m}}{2}(1-x^2) \right)\\
  &=\frac{1}{4}\max_{x\in [-1,1]}\left( ax^2+2 b x + c \right),
\end{align*}
where $a := \delta_{m,k}+\delta_{m,l}-2\sqrt{m}$,  $b := \delta_{m,k}-\delta_{m,l}$, and $c:= \delta_{m,k}+\delta_{m,l}+2\sqrt{m}$.
Taking momentarily for granted that $a < 0$ and that $x_* := -b / a \in [-1,1]$,
we deduce that
\begin{align*}
 \lambda(m,k+l) & \ge
 \frac{1}{4} \left( ax_*^2+2 b x_* + c \right)
 =  \frac{1}{4}\left( - \frac{b^2}{a} + c \right)
 = \frac{1}{4} \frac{b^2 - a c}{-a}\\
 & = \frac{1}{4} \frac{(\delta_{m,k}-\delta_{m,l})^2 + (2\sqrt{m} - \delta_{m,k} -\delta_{m,l}) (2\sqrt{m} + \delta_{m,k}+\delta_{m,l})  }{2\sqrt{m} - \delta_{m,k} - \delta_{m,l}}\\
& = \frac{1}{4} \frac{4m-4\delta_{m,k}\delta_{m,l}}{2\sqrt{m}-\delta_{m,k}-\delta_{m,l}},
\end{align*}
which is the announced lower bound.
It now remains to notice that $a < 0$ and that $-b / a \in [-1,1]$,
but both follow from the general observation that
$\delta_{m,n} < \sqrt{m}$ for $n \ge m > 1$, see Remark \ref{RemGam}. \qed

Before uncovering a family of mutually unbiased equiangular tight frames in the next section, we emphasize here two noteworthy properties relating the vector~$t_\theta$ and the matrix $U_\theta$ that appeared in the above proof.

\begin{prop}
Let $\gamma_{m,k,l}$ be the lower bound for $\lambda(m,k+l)$ from Theorem~\ref{amazing}
and let $\theta \in [0,\pi/2]$ be the angle used in its proof, i.e.,
$$
\gamma_{m,k,l} = \frac{m-\delta_{m,k}\delta_{m,l}}{2\sqrt{m}-\delta_{m,k}-\delta_{m,l}}
\qquad \mbox{and} \qquad
\cos(2\theta) = \frac{\delta_{m,k} - \delta_{m,l}}{2\sqrt{m}-\delta_{m,k}-\delta_{m,l}}.
$$
Then, with $t_\theta \in \mathbb{R}^{k+l},U_\theta \in \mathbb{R}^{m \times (k+l)}$ defined as in \eqref{DeftU} and with $T_\theta:= {\rm diag}[t_\theta]$, one~has
\begin{align}
\label{SP1}
  |U_\theta^\top U_\theta | \, t_\theta  & = \gamma_{m,k,l} \, t_\theta,\\
\label{SP2}
  T_\theta \sgn(U_\theta^\top U_\theta) T_\theta \, U_\theta^\top & = \frac{\gamma_{m,k,l}}{m} \, U_\theta^\top.
\end{align}
\end{prop}

\noindent
{\sc Proof.}
When establishing both \eqref{SP1} and \eqref{SP2}, it will be useful to keep in mind that $\delta_{m,n}$ is tied to $\phi_{m,n} =  \sqrt{(n-m)/(m(n-1))}$ via
$$
\delta_{m,n}
= \frac{m}{n} \bigg( 1 + (n-1) \phi_{m,n} \bigg)
= \frac{m}{n} \bigg( 1 + \frac{n-m}{m} \frac{1}{\phi_{m,n}} \bigg).
$$

Starting with the justification of \eqref{SP1}, we notice that,
since the matrix $V^\top V$ has diagonal entries equal to $1$
and off-diagonal entries equal to $\phi_{m,k}$ in absolute value,
we have
$$
|V^\top V| = (1-\phi_{m,k}) {\rm I}_k + \phi_{m,k} \mathbb{1}_{k,k},
$$
where $\mathbb{1}_{n,n'}$ denotes the $n \times n'$ matrix with all entries equal to $1$.
It follows that
$$
|V^\top V| \mathbb{1}_{k} = (1-\phi_{m,k}) \mathbb{1}_{k} + k \phi_{m,k} \mathbb{1}_{k} = (1+(k-1) \phi_{m,k}) \mathbb{1}_{k}
= \frac{k}{m} \delta_{m,k} \mathbb{1}_{k}.
$$
Likewise, we can obtain
$$
|W^\top W| \mathbb{1}_{l} =  \frac{l}{m} \delta_{m,l} \mathbb{1}_{l}.
$$
Moreover, since the matrices $V^\top W$ and $W^\top V$ have entries all equal to $1/\sqrt{m}$ in absolute value,
we have $| V^\top W | = (1 / \sqrt{m}) \mathbb{1}_{k,l}$
and $| W^\top V | = (1 / \sqrt{m}) \mathbb{1}_{l,k}$,
so that
$$
| V^\top W | \mathbb{1}_{l} = \frac{l}{\sqrt{m}} \mathbb{1}_{k}
\qquad \mbox{and} \qquad
| W^\top V | \mathbb{1}_{k} = \frac{k}{\sqrt{m}} \mathbb{1}_{l}.
$$
Therefore, according to the block-forms of $t_\theta$ and $U_\theta^\top U_\theta$ (see \eqref{DeftU} and \eqref{UTU}),
we observe that
\begin{align}
\nonumber
|U_\theta^\top U_\theta | \, t_\theta & = \begin{bmatrix}
\cos^2 \theta \dfrac{m}{k} \cos \theta \dfrac{1}{\sqrt{k}} \dfrac{k}{m} \delta_{m,k} \mathbb{1}_{k}
+ \cos \theta \sin \theta \dfrac{m}{\sqrt{kl}} \sin \theta \dfrac{1}{\sqrt{l}} \dfrac{l}{\sqrt{m}} \mathbb{1}_{k}
\\
\hline
\cos \theta \sin \theta \dfrac{m}{\sqrt{kl}} \cos \theta \dfrac{1}{\sqrt{k}} \dfrac{k}{\sqrt{m}} \mathbb{1}_{l} + \sin^2 \theta \dfrac{m}{l} \sin \theta \dfrac{1}{\sqrt{l}} \dfrac{l}{m} \delta_{m,l}  \mathbb{1}_{l}
\end{bmatrix}
\\
\label{AlmostThereSP1}
& = \begin{bmatrix}
\cos \theta \dfrac{1}{\sqrt{k}}
\left( \cos^2 \theta \delta_{m,k} + \sin^2 \theta \sqrt{m} \right) \mathbb{1}_{k}\\
\hline
\sin \theta \dfrac{1}{\sqrt{l}}
\left( \cos^2 \theta \sqrt{m} + \sin^2 \theta \delta_{m,l}
\right)
\mathbb{1}_{l}
\end{bmatrix}.
\end{align}
Next, in view of
\begin{align*}
\cos^2 \theta & =
\frac{1+\cos(2 \theta)}{2}
= \frac{\sqrt{m} - \delta_{m,l} }{2 \sqrt{m} - \delta_{m,k} - \delta_{m,l}},\\
\sin^2 \theta & =
\frac{1-\cos(2 \theta)}{2}
= \frac{\sqrt{m} - \delta_{m,k} }{2 \sqrt{m} - \delta_{m,k} - \delta_{m,l}},
\end{align*}
we easily derive that
\begin{equation}
\label{TrigEq}
\cos^2 \theta \delta_{m,k} + \sin^2 \theta \sqrt{m}
= \cos^2 \theta \sqrt{m} + \sin^2 \theta \delta_{m,l}
= \gamma_{m,k,l}.
\end{equation}
When substituting the latter into \eqref{AlmostThereSP1},
the identity \eqref{SP1} immediately follows.

Turning now to the justification of \eqref{SP2},
recalling that
the matrix $V^\top V$ has diagonal entries equal to $1$
and off-diagonal entries equal to $\phi_{m,k}$ in absolute value,
the diagonal entries of the matrix $\sgn(V^\top V)$ are equal to $1$
and its off-diagonal entries are equal to those of $V^\top V$ divided by $\phi_{m,k}$.
In short, we see that $\sgn(V^\top V) = (1-1/\phi_{m,k}){\rm I}_k + (1/\phi_{m,k}) V^\top V$ holds,
and a similar identity holds for $\sgn(W^\top W)$.
Moreover, we also have $\sgn(V^\top W) = \sqrt{m} \, V^\top W$
and $\sgn(W^\top V) = \sqrt{m} \, W^\top V$,
as a consequence of all the entries of $W^\top V$ and $W^\top V$ being equal to $1/\sqrt{m}$ in absolute value.
All in all, according to the block-form~\eqref{UTU} of $U_\theta^\top U_\theta$, we obtain
$$
{\small
\sgn(U_\theta^\top U_\theta)
= \begin{bmatrix}
\left( 1-\dfrac{1}{\phi_{m,k}} \right) {\rm I}_k + \dfrac{1}{\phi_{m,k}} V^\top V & \vline & \sqrt{m} \, V^\top W  \\
\hline
\sqrt{m} \, W^\top  V & \vline & \left( 1-\dfrac{1}{\phi_{m,l}} \right) {\rm I}_l + \dfrac{1}{\phi_{m,l}} W^\top W
\end{bmatrix}.
}
$$
In turn, using the block-form of $T_\theta = {\rm diag}[t_\theta]$,
we derive that $T_\theta \sgn(U_\theta^\top U_\theta)  T_\theta$
takes the form
$$
{\footnotesize
\begin{bmatrix}
\cos^2 \theta \dfrac{1}{k} \left( \left( 1-\dfrac{1}{\phi_{m,k}} \right) {\rm I}_k + \dfrac{1}{\phi_{m,k}} V^\top V \right) & \vline & \cos \theta \sin \theta \dfrac{1}{\sqrt{kl}} \sqrt{m} \, V^\top W  \\
\hline
\cos \theta \sin \theta \dfrac{1}{\sqrt{kl}} \sqrt{m} \, W^\top  V &
\vline & \sin^2 \theta \dfrac{1}{l} \left( \left( 1-\dfrac{1}{\phi_{m,l}} \right) {\rm I}_l + \dfrac{1}{\phi_{m,l}} W^\top W  \right)
\end{bmatrix}.
}
$$
Multiplying on the right by the transpose of $U_\theta= {\small \begin{bmatrix} \; \cos\theta \sqrt{\dfrac{m}{k}} V & \vline & \sin \theta \sqrt{\dfrac{m}{l}} W  \; \end{bmatrix}}$
and making use of the facts that $V V^\top = (k/m) {\rm I}_m$
and $W W^\top = (l/m) {\rm I}_m$,
the matrix $T_\theta  \sgn(U_\theta^\top U_\theta) T_\theta \, U_\theta^\top$ becomes
\begin{align*}
 & {\footnotesize \begin{bmatrix}
\cos^2 \theta \dfrac{1}{k} \cos \theta \sqrt{\dfrac{m}{k}}
\left( \left( 1-\dfrac{1}{\phi_{m,k}} \right) + \dfrac{k}{m} \dfrac{1}{\phi_{m,k}} \right) V^\top
+ \cos \theta \sin \theta \dfrac{1}{\sqrt{k l}} \sqrt{m} \sin \theta \sqrt{\dfrac{m}{l}} \dfrac{l}{m} V^\top\\
\hline
\cos \theta \sin \theta \dfrac{1}{\sqrt{kl}} \sqrt{m} \cos\theta \sqrt{\dfrac{m}{k}} \dfrac{k}{m} W^\top
+ \sin^2 \theta \dfrac{1}{l} \sin \theta \sqrt{\dfrac{m}{l}}
\left( \left( 1-\dfrac{1}{\phi_{m,l}} \right)  + \dfrac{l}{m} \dfrac{1}{\phi_{m,l}} \right) W^\top
\end{bmatrix} }\\
& = \begin{bmatrix}
\cos \theta \sqrt{\dfrac{m}{k}} \left( \dfrac{\cos^2 \theta}{k} \left( 1 + \dfrac{k-m}{m} \dfrac{1}{\phi_{m,k}} \right)
+ \sin^2 \theta \dfrac{1}{\sqrt{m}}
\right) V^\top \\
\hline
\sin \theta \sqrt{\dfrac{m}{l}} \left(
\cos^2 \theta \dfrac{1}{\sqrt{m}} + \dfrac{\sin^2 \theta}{l} \left( 1 + \dfrac{l-m}{m} \dfrac{1}{\phi_{m,l}}\right)
\right) W^\top
\end{bmatrix}\\
& = \begin{bmatrix}
\cos \theta \sqrt{\dfrac{m}{k}} \left( \dfrac{\cos^2 \theta}{m}
 \delta_{m,k}
+  \dfrac{\sin^2 \theta}{\sqrt{m}}
\right) V^\top \\
\hline
\sin \theta \sqrt{\dfrac{m}{l}} \left(
\dfrac{\cos^2 \theta}{\sqrt{m}} + \dfrac{\sin^2 \theta}{m} \delta_{m,l}
\right) W^\top
\end{bmatrix}.
\end{align*}
Similarly to \eqref{SP1},
the identity \eqref{SP2} now simply
follows by exploiting \eqref{TrigEq} again.
\qed

\section{Construction of Mutually Unbiased Equiangular Tight Frames}

To apply the result of Theorem~\ref{amazing} in practical situations,
we evidently need to uncover specific integers $k$, $l$, and $m$ allowing mutually unbiased equiangular tight frames to exist.
As a simple example, one can take $k=l=m$ and consider $(v_1,\ldots,v_k)$
to be the canonical basis for $\mathbb{R}^m$
and $(w_1,\ldots,w_l)$ to be the columns of an $m \times m$ Hadamard matrix
--- recall that $m \times m$ Hadamard matrices are conjectured to exist when and only when $m$ is a multiple of $4$ (the `only when' part being acquired, of course).
This would yield the lower bound $\lambda(m) \ge (1+\sqrt{m})/2$, $m \in 4 \mathbb{N}$,
which is inferior to the lower bounds reported
in \cite{FS} for $m=4$ and $m=8$.
As a slightly more elaborate example,
one can take $k=m$ and $(v_1,\ldots,v_k)$ to be the canonical basis of $\mathbb{R}^m$,
together with $l > m$ and $(w_1,\ldots,w_l)$ to be a real equiangular tight frame for $\mathbb{R}^m$ that is \textit{flat},
in the sense that every entry of each vector $w_j$ is either $1/\sqrt{m}$ or $-1/\sqrt{m}$.
Real flat equiangular tight frames are equivalent to binary codes achieving equality in the Grey--Rankin bound
and infinite families  are known (see \cite{JasperMF14,FickusJMP21}).
This would yield the lower bound
$\lambda(m,m+l) \geq (m-\gamma_{m,l})/(2\sqrt{m}-1-\gamma_{m,l})$.
With $m=6$ and $l=16$,
this provides the lower bound $\lambda(6)\gtrsim 2.2741$,
which is superior to the lower bounds reported in~\cite{FS}
but inferior to the numerical evaluation $\lambda(6) \approx 2.2857$
performed by B.~L.~Chalmers and corroborated by our own computations.
In order to apply Theorem~\ref{amazing} more effectively,
we need further examples of mutually unbiased equiangular tight frames.
To this end, we now relate such frames to a type of generalized Hadamard matrices.

\begin{prop}
\label{muetfCharacterization}
Given integers $k,l \geq m > 1$,
there are mutually unbiased equiangular tight frames $(v_1,\dots,v_k)$ and $(w_1,\dots,w_l)$ for~$\mathbb{R}^m$ if and only if there is a $k\times l$ matrix $X$ with the following five properties:
\begin{enumerate}[(i)]
\item \label{P1}
$X_{ij}\in\{-1,+1\}\,$ for all $i\in\{1,\dots,k\}$ and $j\in \{1,\dots,l\}$;
\item \label{P2}
$XX^\top X=aX\,$ for some $a\in\mathbb{R}$;
\item \label{P3}
$X$ has equiangular rows, i.e., $|X X^\top|_{i,i'}$ is constant over all $i \neq i'$;
\item \label{P4}
$X$ has equiangular columns, i.e.,
$|X^\top X|_{j,j'}$ is constant over all $j \neq j'$;
\item \label{P5}
$X$ has rank $m$.
\end{enumerate}
When this occurs, the following three quantities are necessarily integers:
\begin{equation}
\label{integrality}
\frac{k l}{m},
\qquad
k\,\sqrt{\frac{l-m}{m(l-1)}},
\qquad
l\,\sqrt{\frac{k-m}{m(k-1)}}.
\end{equation}
\end{prop}

\noindent
{\sc Proof.}
Firstly, let us assume that there are mutually unbiased equiangular tight frames $(v_1,\dots,v_k)$ and $(w_1,\dots,w_l)$ for~$\mathbb{R}^m$.
With $V \in \mathbb{R}^{m \times k}$ and $W \in \mathbb{R}^{m \times l}$ denoting the matrices with columns $v_1,\dots,v_k$ and $w_1,\dots,w_l$, respectively,
we set
$$
X=\sqrt{m}\,V^\top W\in\mathbb{R}^{k\times l}.
$$
By Lemma~\ref{beta}, we have $|V^\top W|_{i,j} = |\langle v_i, w_j \rangle | = 1/\sqrt{m}$ for all $i\in\{1,\dots,k\}$ and $j\in \{1,\dots,l\}$,
so Property \eqref{P1} is immediate.
In view of $VV^\top = (k/m)\,{\rm I}_m$ and of $WW^\top=(l/m)\,{\rm I}_m$,
it is also straightforward to see that
\begin{equation}
\label{XX'}
XX^\top = l \, V^\top V
\qquad \mbox{and} \qquad
X^\top X =k \, W^\top W.
\end{equation}
From here, using the fact that $VV^\top = (k/m)\,{\rm I}_m$ one more time,
we obtain that $XX^\top X = (l \, V^\top V) (\sqrt{m}\, V^\top W)= (kl/m)\sqrt{m}\, V^\top W$,
i.e., $XX^\top X  = a X$ with $a = kl/m$,
so Property \eqref{P2} is satisfied.
Properties \eqref{P3} and \eqref{P4}, too, are consequences of~\eqref{XX'},
since e.g. the off-diagonal entries of $X X^\top$ are constant in absolute value because those of $V^\top V$ are.
Finally, Property \eqref{P5} is also implied  by~\eqref{XX'} via~${\rm rank}(X) = {\rm rank}(X X^\top) = {\rm rank}(V^\top V) = {\rm rank}(V V^\top) = {\rm rank}({\rm I}_m)~=~m$.

Conversely, let us assume that Properties \eqref{P1}--\eqref{P2} are fulfilled by some matrix $X \in \mathbb{R}^{k \times l}$.
Consider the singular value decomposition of this matrix written as $X = P \Sigma Q^\top$, where the diagonal matrix $\Sigma \in \mathbb{R}^{m \times m}$ has positive entries (by \eqref{P5})
and where the matrices $P \in \mathbb{R}^{k \times m}$ and $Q \in \mathbb{R}^{l \times m}$ have orthonormal columns,
i.e., $P^\top P = {\rm I}_m$ and $Q^\top Q = {\rm I}_m$.
Property \eqref{P2} easily yields $\Sigma^3 = a \, \Sigma$ and hence $\Sigma = \sqrt{a} \, {\rm I}_m$.
Then, looking at the squared Frobenius norm of $X = \sqrt{a} \, P Q^\top$,
we derive from \eqref{P1} that $k l = a m$, i.e., that $a = kl/m$.
We now set
$$
V = \sqrt{\frac{k}{m}}P^\top \in \mathbb{R}^{m \times k}
\qquad \mbox{and} \qquad
W = \sqrt{\frac{l}{m}} Q^\top \in \mathbb{R}^{m \times l}
$$
and we claim that the columns $v_1,\ldots,v_k$ of $V$
and $w_1,\ldots,w_l$ of $W$ are mutually unbiased equiangular tight frames for $\mathbb{R}^m$.
Indeed, using $V^\top V = (k/m) P P^\top$
and $X X^\top = a P P^\top$,
we see that $V^\top V = (1/l) X X^\top$,
so that the equiangularity of the system $(v_1,\ldots,v_k)$ is clear from \eqref{P3}.
Note that each $v_i$ is a unit vector,
since $\|v_i\|_2^2 = (V^\top V)_{i,i} = (1/l) (XX^\top)_{i,i} = (1/l) \sum_{j=1}^l X_{i,j}^2 = 1$ by \eqref{P1}.
The fact that these vectors form a tight frame is  seen from $V V^\top = (k/m) \, P^\top P = (k/m) \, {\rm I}_m$.
Similar arguments (using \eqref{P4}) would reveal that the system $(w_1,\ldots,w_l)$ is also an equiangular tight frame.
At last, to see that these systems are mutually unbiased,
it suffices to notice that $V^\top W = (\sqrt{kl}/m) \, P Q^\top = (1/\sqrt{m}) \, X$
and to invoke \eqref{P1} once again.

It finally remains to establish that the three quantities in \eqref{integrality} are integers.
For the first one,
we have seen  (in the proofs of both implications)
that $a = kl/m$
and \eqref{P1}-\eqref{P2} show that $a$ is an integer:
any entry of $XX^\top X=a X$ is on the one hand an integer and on the other hand equal to $\pm a$.
For the third one, say,
looking e.g. at~\eqref{XX'},
any off-diagonal entry of $X X^\top = l \, V^\top V$
is on the one hand an integer and on the other hand equal to $l$ times the common absolute inner product in a $k$-vector equiangular tight frame for $\mathbb{R}^m$,
i.e., to $l \sqrt{(k-m)/(m(k-1))}$.
\qed
\vskip .1in

Although conditions \eqref{P1}--\eqref{P5} are restrictive,
there are matrices $X$ satisfying them with $ m < k < l$.
For instance, the $6\times 10$ matrix
\begin{equation*}
\setlength{\arraycolsep}{1pt}
X=\begin{small}\left[\begin{array}{rrrrrrrrrr}
 1& 1& 1& 1& 1& 1& 1& 1& 1& 1\\
 1& 1&-1& 1&-1&-1& 1&-1&-1&-1\\
 1&-1& 1&-1& 1&-1&-1& 1&-1&-1\\
-1& 1& 1&-1&-1& 1&-1&-1& 1&-1\\
-1&-1&-1& 1& 1& 1&-1&-1&-1& 1\\
-1&-1&-1&-1&-1&-1& 1& 1& 1& 1
\end{array}\right]\end{small}
\end{equation*}
is one such matrix\footnote{As pointed out to us by Josiah~Park, this same $6\times 10$ matrix appeared in a recent investigation of spherical half-designs (see~\cite{HW}).}:
it has $\pm 1$ entries,
the identity $X X^\top X = a X$ is easily verified (at least computationally),
and it was already observed in~\cite{FickusMJ16}
that both its rows and its columns form equiangular tight frames for their $5$-dimensional spans.
Therefore, since $X$ fulfills the conditions of Proposition~\ref{muetfCharacterization} with $m=5$, $k=6$, and $l=10$,
we are guaranteed the existence of  mutually unbiased equiangular tight frames $(v_1,\dots,v_6)$ and $(w_1,\dots,w_{10})$ for $\mathbb{R}^5$.
Remarkably, this example is but the first member of the infinite family presented below.

\begin{thm}
\label{muetfFamily}
For any integer $s \ge 2$, there are mutually unbiased equiangular tight frames $(v_1,\dots,v_k)$ and $(w_1,\dots,w_l)$ for~$\mathbb{R}^m$, where
\begin{equation*}
k = 2^{s-1}(2^s-1),
\qquad
l = 2^{s-1}(2^s+1),
\qquad
m = \frac{2^{2s}-1}{3}.
\end{equation*}
\end{thm}

\noindent
{\sc Proof.}
For any such $s$, $k$, $l$ and $m$,
the requisite matrix $X$ of Proposition~\ref{muetfCharacterization} is produced in the recent paper~\cite{FickusIJK21},
albeit nonobviously so.
In brief, let $Q$ and $B$ be the canonical hyperbolic-quadratic and symplectic forms on the binary vector space $\mathbb{F}_2^{2s}$, respectively:
\begin{align*}
Q(x)=Q(x_1,\dotsc,x_{2s})
&:=\sum_{r=1}^s x_{2r-1}x_{2r},\\
B(x,y)=B((x_1,\dotsc,x_{2s}),(y_1,\dotsc,y_{2s}))
&:=\sum_{r=1}^s (x_{2r-1}y_{2r}+x_{2r}y_{2r-1}).
\end{align*}
Let $\Gamma$ be the corresponding character table of $\mathbb{F}_2^{2s}$,
defined by $\Gamma(x,y)=(-1)^{B(x,y)}$ for all $x,y\in\mathbb{F}_2^{2s}$.
Any submatrix of $\Gamma$ obviously satisfies \eqref{P1} from Proposition~\ref{muetfCharacterization}.
Let $X$ be the specific submatrix of $\Gamma$ whose rows and columns are indexed by
$\{x\in\mathbb{F}_2^{2s}: Q(x)=1\}$ and $\{x\in\mathbb{F}_2^{2s}: Q(x)=0\}$, respectively.
By Lemma~4.2 of~\cite{FickusIJK21}, these two subsets of $\mathbb{F}_2^{2s}$ are \textit{difference sets} for $\mathbb{F}_2^{2s}$ of cardinality $k$ and $l$, respectively.
As detailed in~\cite{FickusIJK21}, this means that the rows and columns of $X$ are equiangular, namely that \eqref{P3} and \eqref{P4} hold.
Theorem~4.4 of~\cite{FickusIJK21} moreover gives that these two difference sets are \textit{paired}, meaning that the columns of $X$ form a tight frame for their span,
so that \eqref{P2} holds.
Theorem~3.3 of~\cite{FickusIJK21} then implies that the rank of $X$ is indeed $m$, so that \eqref{P5} holds.
\qed
\vskip .1in

We close this section by highlighting that real mutually unbiased equiangular tight frames are rare objects.
Precisely, we have obtained rather stringent necessary conditions for their existence (not included here because too detached from our main focus).
For instance, these conditions imply that mutually unbiased $k$-vector and $l$-vector equiangular tight frames for $\mathbb{R}^m$ can only exist for at most thirteen triples of integers $(m,k,l)$ with $l > k > m+1$ when $m \le 1000$,
and that they cannot exist when $l=k>m$, in contrast with the complex setting.

\section{Epilogue: the fifth maximal projection constant}

By combining the main results derived in the two previous sections,
namely Theorems~\ref{amazing} and~\ref{muetfFamily},
and after some tedious algebraic manipulation,
we can state that the maximal relative projection constant
at any $m$ of the form $m = (2^{2s}-1)/3$ for some integer $s \ge 2$
is bounded below as
\begin{equation}
\label{LBFam}
\lambda ( m ,4^s) \ge
\frac{2^{2s}-1}{2^{3s} - 3\, 2^{s-1} + 1}
\left( \frac{2^{2s-1} + 2^s -1}{3} + 2^{s-1} \sqrt{m} \right).
\end{equation}
If this was to be an equality,
then the vector $t_\theta \in \mathbb{R}_+^{n}$, $n = 4^s$,
and the matrix $U_\theta \in \mathbb{R}^{m \times n}$, $m = (2^{2s}-1)/3$,
appearing in the proof of Theorem \ref{amazing}
should be maximizers of the expression for $\lambda(m,n)$ from Theorem \ref{lammbda}.
For genuine maximizers $\bar{t} \in \mathbb{R}_+^{n}$
and $\bar{U} \in \mathbb{R}^{m \times n}$, we emphasize the following two necessary conditions:

\begin{enumerate}[(a)]

\item \label{NecCond1}
$\bar{t}$ is a maximizer of $\sum_{i,j} t_i t_j |\bar{U}^\top \bar{U} | _{i,j}$ subject to $\|t\|_2 = 1$,
so is characterized by the fact that $\bar{t}$
is an eigenvector (in fact, the leading eigenvector) of $|U^\top U|$ --- this is indeed satisfied by $t_\theta$ and $U_\theta$, according to \eqref{SP1};

\item \label{NecCond2}
$\bar{U}$ is a maximizer of $\sum_{i,j} \bar{t}_i \bar{t}_j \sgn(\bar{U}^\top \bar{U})_{i,j} (U^\top U) _{i,j} = \textrm{tr} (\bar{T} \sgn(\bar{U}^\top \bar{U}) \, \bar{T} U^\top U)$,
$\bar{T} := {\rm diag}[\bar{t}]$,
 subject to $U U^\top = {\rm I}_m$,
so is characterized by the fact that the rows of $\bar{U}$ are eigenvectors corresponding to the $m$ largest eigenvalues of
$\bar{T} \sgn(\bar{U}^\top \bar{U}) \bar{T}$ ---
 this is indeed satisfied by $t_\theta$ and $U_\theta$, according to \eqref{SP2}.

\end{enumerate}

\begin{rem}
\label{RkCombine2NedCond}
The necessary conditions \eqref{NecCond1}-\eqref{NecCond2} combine to show that the genuine maximizers $\bar{t}$ and $\bar{U}$ obey the noteworthy relation
$$
\big( \bar{U}^\top \bar{D} \bar{U} \big)_{i,i}
= \lambda(m,n) \, \bar{t}_i^2
\qquad \mbox{for all } i \in \{1,\ldots,n\},
$$
where $\bar{D} = {\rm diag}[\bar{\mu}_1,\ldots,\bar{\mu}_m]$ is the diagonal matrix with the $m$ leading eignevalues $\bar{\mu}_1 \ge \cdots \ge \bar{\mu}_m$ of $\bar{T} \sgn(\bar{U}^\top \bar{U}) \bar{T}$ on its diagonal.
Indeed, by \eqref{NecCond1}, we have
\begin{align}
\nonumber
\lambda(m,n) \, \bar{t}_i^2
& = \bar{t}_i \sum_{j=1}^n |\bar{U}^\top \bar{U}|_{i,j} \bar{t}_j
= \sum_{j=1}^n (\bar{U}^\top \bar{U})_{i,j}
(\bar{T} \sgn(\bar{U}^\top \bar{U}) \bar{T})_{i,j}\\
\label{Forti2=}
& = \big( (\bar{U}^\top \bar{U}) (\bar{T} \sgn(\bar{U}^\top \bar{U}) \bar{T}) \big)_{i,i}.
\end{align}
Now, by \eqref{NecCond2}, we have
$\bar{T} \sgn(\bar{U}^\top \bar{U} \bar{T}) \bar{U}^\top = \bar{U}^\top \bar{D}$,
or $\bar{U} \bar{T} \sgn(\bar{U}^\top \bar{U}) \bar{T}  =  \bar{D} \bar{U}$ by taking the transpose.
Making use of the latter in \eqref{Forti2=} gives the expected relation.
\end{rem}

The observation that $t_\theta$ and $U_\theta$ do satisfy conditions \eqref{NecCond1}-\eqref{NecCond2}
supports the belief that \eqref{LBFam} could be an equality.
To the question of whether
the right-hand side of \eqref{LBFam} also coincides with the value of the maximal absolute projection constant $\lambda ( m )$, $m = (2^{2s}-1)/3$,
the answer is in general {\em no}.
Indeed, for $s=3$,
hence for $m=21$, $k=28$, and $l=36$,
we have $\gamma_{21,28,36} \approx 3.9397$, while a real equiangular tight frame for $\mathbb{R}^{21}$ made of $126$ vectors is known to exist (see e.g. \cite{FickusM16}), so Theorem~\ref{ETFs} yields $\lambda(21) \ge \lambda(21,126) \gtrsim 4.3333$.
However, for $s=2$,
hence for $m=5$, $k=6$, and $l=10$,
there are convincing reasons to believe that $\gamma_{5,6,10} \approx 2.06919$ coincide with the value of $\lambda(5)$.
These reasons are the extensive numerical investigations carried out B. L. Chalmers,
as well as our own computations
(some of which can be found in a {\sc matlab} reproducible available on the authors' webpages).
All these clues prompt us to conclude with the following assertion.

\begin{thm}[and Conjecture]
The fifth absolute projection constant satisfies
\begin{equation*}
 \lambda(5) \ge  \lambda(5,16)\geq \frac{5}{59}(11+6\sqrt{5}) \approx 2.06919,
\end{equation*}
and it is expected that the latter is indeed the true value of $\lambda(5)$.
\end{thm}

\vskip .5in

\section*{Appendix}

As bonus material,
we present here a new proof of Theorem \ref{ETFs}
as a immediate consequence of the technical result below coupled with Theorem \ref{lammbda}.

\begin{prop}
\label{PropEquiv}
For integers $n \ge m >1$,
one has
\begin{align}
\nonumber
\max  \bigg\lbrace \sum_{i,j=1}^n &  t_it_j|U^\top U|_{ij}:t\in\mathbb{R}^n,\;\|t\|_2=1,U\in \mathbb{R}^{m\times n},\; UU^\top={\rm I}_m \bigg\rbrace \\
\label{UB_proved}
& \le
\frac{m}{n} \left( 1 + \sqrt{\frac{(n-1)(n-m)}{m}} \right),
\end{align}
with equality if and only if there exists a matrix $U \in \mathbb{R}^{m \times n}$ with $U U^\top  = {\rm I}_m$,
$(U^\top U)_{i,i} = m/n$
for all $i  \in \{1,\ldots, n\}$,
and  $|U^\top U|_{i,j} = \sqrt{(n-m)m/(n-1)}/n$ for all $i \not= j \in \{1,\ldots, n\}$.
\end{prop}

\noindent
{\sc Proof.}
For $t \in \mathbb{R}^n$ satisfying $\|t\|_2 = 1$ and $U \in \mathbb{R}^{m \times n}$ satisfying $U U^\top = {\rm I}_m$,
we  use the nonnegativity of $(U^\top U)_{i,i}$ (as the inner product of the $i$th column of $U$ with itself) and Cauchy--Schwarz inequality
to write
\begin{align*}
\Sigma & :=
\sum_{i,j=1}^n t_i t_j  |U^\top U|_{i,j}
 =  \sum_{i=1}^n t_i^2 |U^\top U|_{i,i}  + \sum_{\substack{i,j=1 \\ i \not= j}}^n t_i t_j |U^\top U|_{i,j} \\
& \le \sum_{i=1}^n t_i^2 (U^\top U)_{i,i}  +  \sqrt{ \sum_{\substack{i,j=1 \\ i \not= j}}^n t_i^2 t_j^2  } \sqrt{ \sum_{\substack{i,j=1 \\ i \not= j}}^n  (U^\top U)_{i,j}^2 }  \\
& =  \sum_{i=1}^n t_i^2 (U^\top U)_{i,i} + \sqrt{ \sum_{i,j=1}^n t_i^2 t_j^2 - \sum_{i=1}^n t_i^4  } \sqrt{ \sum_{i,j=1}^n  (U^\top U)_{i,j}^2 - \sum_{i=1}^n  (U^\top U)_{i,i}^2 }\\
& = \sum_{i=1}^n \alpha_i \beta_i + \sqrt{ A - \sum_{i=1}^n \alpha_i^2} \sqrt{ B - \sum_{i=1}^n \beta_i^2 },
\end{align*}
where we have set $\alpha_i = t_i^2$,
$\beta_i = (U^\top U)_{i,i}$,
$A = \big(\sum_i t_i^2 \big)  \big(\sum_j t_j^2 \big) = \|t\|_2^4 = 1$,
and $B = \sum_{i,j} (U^\top U)_{i,j}^2 = \|U^\top U\|_F^2 = \textrm{tr}(U^\top U U^\top U) = \textrm{tr}(U U^\top U U^\top) = m $.
Setting also $a = \|t\|_2^2 =1$,  $b = \textrm{tr}(U^\top U) = \textrm{tr}(U U^\top) = m $,
as well as
$$
x_i := \frac{\alpha_i - a/n}{\sqrt{A - a^2/n}}
\qquad \mbox{and} \qquad
y_i := \frac{\beta_i - b/n}{\sqrt{B - b^2/n}},
$$
we notice that $\sum_{i=1}^n x_i =0$ and $\sum_{i=1}^n y_i =0$.
We exploit these identities a few times to derive
\begin{align*}
\Sigma
& \le \sum_{i=1}^n \bigg( \frac{a}{n} + \sqrt{A - \frac{a^2}{n}} x_i \bigg) \bigg( \frac{b}{n} + \sqrt{B - \frac{b^2}{n}} y_i \bigg)\\
& + \sqrt{A - \sum_{i=1}^n \bigg( \frac{a}{n} + \sqrt{A - \frac{a^2}{n}} x_i \bigg)^2}
\sqrt{ B - \sum_{i=1}^n \bigg( \frac{b}{n} + \sqrt{B - \frac{b^2}{n}} y_i \bigg)^2}\\
& = \frac{ab}{n} + \sqrt{A - \frac{a^2}{n}} \sqrt{B - \frac{b^2}{n}}  \sum_{i=1}^n x_i y_i\\
& + \sqrt{A - \frac{a^2}{n} - \Big( A - \frac{a^2}{n} \Big) \sum_{i=1}^n x_i ^2}
+ \sqrt{B - \frac{b^2}{n} - \Big( B - \frac{b^2}{n} \Big) \sum_{i=1}^n y_i ^2} \\
& = \frac{ab}{n} + \sqrt{A - \frac{a^2}{n}} \sqrt{B - \frac{b^2}{n}}
\left[ \sum_{i=1}^n x_i y_i  + \sqrt{1-\sum_{i=1}^n x_i^2} \sqrt{1-\sum_{i=1}^n y_i^2} \right].
\end{align*}
The latter term in square brackets is nothing but the inner product of the unit vectors $\tilde{x} := \begin{bmatrix} x , \sqrt{1-\|x\|_2^2} \end{bmatrix}$ and $\tilde{y} := \begin{bmatrix} y , \sqrt{1-\|y\|_2^2} \end{bmatrix}$, so it is bounded by one.
Thus, keeping the values of $a=1$, $b=m$, $A=1$, and $B=m$ in mind,
we arrive at
$$
\sum_{i,j=1}^n t_i t_j  |U^\top U|_{i,j}
\le \frac{m}{n} + \sqrt{1-\frac{1}{n}} \sqrt{m-\frac{m^2}{n}}.
$$
Taking the supremum over $t$ and $U$ leads to the desired inequality \eqref{UB_proved} after some algebraic manipulation.
This inequality turns into an equality if the matrix $U \in \mathbb{R}^{m \times n}$ with $U U^\top  = {\rm I}_m$ satisfies
$(U^\top U)_{i,i} = m/n$
for all $i  \in \{1,\ldots, n\}$
and  $|U^\top U|_{i,j} = \sqrt{(n-m)m/(n-1)}/n$ for all $i \not= j \in \{1,\ldots, n\}$,
simply by choosing $t \in \mathbb{R}^{n}$ with entries $t_i = 1/\sqrt{n}$
for all $i \in \{ 1, \ldots,n\}$.

Conversely,  let us assume that \eqref{UB_proved} is an equality.
Our goal is now to prove that
$(U^\top U)_{i,i} = m/n$ for all $i  \in \{1,\ldots, n\}$
and  $|U^\top U|_{i,j} = \sqrt{(n-m)m/(n-1)}/n$ for all $i \not= j \in \{1,\ldots, n\}$,
where $U \in \mathbb{R}^{m \times n}$ satisfying $U U^\top = {\rm I}_m$
achieves the maximum,
together with $t \in \mathbb{R}^n$ satisfying $\|t\|_2 = 1$.
We start by taking into account that equality must hold throughout the first part of the argument.
Equality in Cauchy--Schwarz inequality implies the existence of $c \in \mathbb{R}$ such that
$$
t_i t_j   = c \, |U^\top U|_{i,j}
\qquad \mbox{for all  }i \not= j \in \{1,\ldots, n\}
$$
and equality in $\langle \tilde{x},\tilde{y} \rangle \le 1$ yields $x=y$,
i.e.,
\begin{equation}
\label{EqIn...}
(U^\top U)_{i,i} -  \frac{m}{n} = \frac{\sqrt{m - m^2/n}}{\sqrt{1-1/n}} \bigg( t_i^2 - \frac{1}{n} \bigg)
 \qquad \mbox{for all  }i  \in \{1,\ldots, n\}.
\end{equation}
Since the matrix $T \sgn(U^\top U) T$ has diagonal entries $(T \sgn(U^\top U) T)_{i,i} = t_i^2$
and off-diagonal entries
$$
(T \sgn(U^\top U) T)_{i,j} = t_i t_j \sgn(U^\top U)_{i,j}
= c \, |U^\top U|_{i,j} \sgn(U^\top U)_{i,j}
= c \, (U^\top U)_{i,j},
$$
the necessary condition \eqref{NecCond2},
written for all $i \in \{1,\ldots,n\}$ and $h \in \{1,\ldots,m\}$ as
$$
\sum_{j=1}^n (T \sgn(U^\top U) T)_{i,j} U^\top_{j,h} = \mu_h  U^\top_{i,h},
$$
where $\mu_1 \ge \cdots \ge \mu_m$ are the $m$ leading eigenvalues of $T (\sgn(U^\top U) T$, becomes
$$
t_i^2 U^\top_{i,h} + \sum_{\substack{j=1\\j\not= i}}^n c \, (U^\top U)_{i,j} U^\top_{j,h}
= \mu_h  U^\top_{i,h}.
$$
In other words, for all $i \in \{1,\ldots,n\}$ and $h \in \{1,\ldots,m\}$,
we have
$$
t_i^2 U^\top_{i,h} + c \, (U^\top U U^\top)_{i,h} - c \, (U^\top U)_{i,i} U^\top_{i,h} = \mu_h  U^\top_{i,h},
$$
or equivalently, in view of $U U^\top = {\rm I}_m$,
\begin{equation}
\label{SomeEq}
\left( t_i^2 + c - c \, (U^\top U)_{i,i} \right) U^\top_{i,h} = \mu_h  U^\top_{i,h}.
\end{equation}
This actually shows that $\mu_h$ is independent of $h \in \{1,\ldots,m\}$
and --- thanks to the alternate expression $\lambda(m,n) = \mu_1+\ldots+\mu_m$ (see e.g. \cite[Theorem 1]{FS}) ---
one must have $\mu_h = \lambda(m,n)/m$.
Now \eqref{SomeEq} reduces (say, by multiplying by $U^\top_{i,h}$, summing over $h$, and simplifying) to
$t_i^2 + c - c \, (U^\top U)_{i,i} = \lambda(m,n)/m$.
Summing  over $i \in \{1,\ldots,n\}$) yields
$$
1 + c \, (n-m) = \frac{n}{m} \lambda(m,n) = 1 + \sqrt{\frac{(n-1)(n-m)}{m}},
$$
which shows that
$$
c= \sqrt{\frac{n-1}{m(n-m)}}.
$$
Invoking Remark \ref{RkCombine2NedCond}, we notice that
$(U^\top U)_{i,i}  = m t_i^2$ for all $i \in \{1,\ldots,n\}$,
and therefore \eqref{EqIn...} becomes $m(t_i^2 - 1/n) = \sqrt{m(n-m)/(n-1)}(t_i^2 - 1/n)$.
Given that $m~\not=~\sqrt{m(n-m)/(n-1)}$ when $m>1$,
we consequently obtain $t_i^2 = 1/n$ for all $i~\in~\{1,\ldots,n\}$.
In turn,
we deduce from $(U^\top U)_{i,i} = m t_i^2$ that $(U^\top U)_{i,i} = m/n$ for all $i \in \{1,\ldots,n\}$
and from $ c |U^\top U|_{i,j} = t_i t_j$ that $|U^\top U|_{i,j} =\sqrt{m(n-m)/(n-1)}/n$ for all $i \not= j \in \{1,\ldots,n\}$.
The proof is now complete.
\qed




\begin{thebibliography}{00}



\bibitem{B} G. Basso, {\it Computation of maximal projection constants,} J. Funct. Anal. 277/10 (2019), 3560--3585.

\bibitem{B2} G. Basso, {\it Almost minimal orthogonal projections,} Isr. J. Math. 243 (2021),  355--376.

\bibitem{CpGaG}
F.  Caro Perez, V. Gonzalez Avella,  D. Goyeneche, {\it Mutually unbiased frames,}
arXiv preprint arXiv:2110.08293 (2021).

\bibitem{CLM} A. Castejon, G. Lewicki, M. Martin, {\it Some results on absolute projection constant,} Numer. Func. Anal. Optim.  40/1 (2019), 34--51.

\bibitem{CLe} B. L. Chalmers, G. Lewicki, \textit{Three-dimensional subspace of $l_{\infty}^{(5)}$ with maximal projection constant,} J. Funct. Anal. 257/2 (2009), 553--592.

\bibitem{CL} B. L. Chalmers, G. Lewicki, {\it A proof of the Gr\"unbaum conjecture,} Studia Math. 200 (2010), 103--129.

\bibitem{FickusIJK21} M.~Fickus, J.~W.~Iverson, J.~Jasper, E.~J.~King, {\it Grassmannian codes from paired difference sets,} Des.\ Codes Cryptogr.\ 89 (2021) 2553--2576.

\bibitem{FickusJMP21} M.~Fickus, J.~Jasper, D.~G.~Mixon, J.~D.~Peterson, {\it Hadamard equiangular tight frames,} Appl.\ Comput.\ Harmon.\ Anal.\ 50 (2021) 281--302.

\bibitem{FM} M. Fickus,  B.  R.  Mayo, {\it Mutually unbiased equiangular tight frames,}
IEEE Trans.\ Inform.\ Theory 67/3 (2020), 1656--1667.

\bibitem{FickusM16} M.~Fickus, D.~G.~Mixon, {\it Tables of the existence of equiangular tight frames,} arXiv:1504.00253 (2016).

\bibitem{FickusMJ16} M.~Fickus, D.~G.~Mixon, J.~Jasper, {\it Equiangular tight frames from hyperovals,} IEEE Trans.\ Inform.\ Theory 62/9 (2016) 5225--5236.

\bibitem{FR} S. Foucart, H. Rauhut, {\it A Mathematical Introduction to Compressive Sensing,} Birkh\"auser, 2013.

\bibitem{FS}S. Foucart, L. Skrzypek, {\it On maximal relative projection constants,} J. Math. Anal. Appl. 447/1 (2017), 309--328.

\bibitem{G} B. Gr\"unbaum, {\it Projection constants,} Trans. Amer. Math. Soc. 95 (1960), 451--465.

\bibitem{HW} D. Hughes, S. Waldron,
{\it Spherical half-designs of high order},
Involve, a Journal of Mathematics 13/2 (2020), 193--203.

\bibitem{JasperMF14} J.~Jasper, D.~G.~Mixon, M.~Fickus, {\it Kirkman equiangular tight frames and codes,} IEEE Trans.\ Inform.\ Theory\ 60/1 (2014) 170--181.

\bibitem{KS} I. M. Kadec, M. G. Snobar, {\it Certain functionals on the Minkowski compactum,} Math. Notes 10 (1971), 694--696 (English transl.).

\bibitem{K} H. K\"{o}nig, {\it Spaces with large projection constants,} Isr. J. Math. 50/3 (1985), 181--188.

\bibitem{KT} H. K\"{o}nig,  N. Tomczak-Jaegermann,
{\it Norms of minimal projections,}
J. Funct. Anal. 119/2 (1994),  253--280.

\bibitem{KLL} H. K\"{o}nig, D. Lewis, P.-K. Lin,
{\it  Finite dimensional projection constants,}
Studia Mathematica 75/3 (1983), 341--358.

\bibitem{S} F. Sokolowski, {\it Minimal projections onto subspaces with codimension 2,} Numer. Funct. Anal. Optim. 38/8 (2017), 1045--1059.


\bibitem{W} P. Wojtaszczyk, {\it Banach Spaces for Analysts,} Cambridge University Press, Cambridge, 1991.


\end{thebibliography}


\end{document}